\documentclass[a4paper,11.5pt,reqno]{amsart}

\usepackage{graphicx}
\usepackage[utf8]{inputenc}
\usepackage{amsmath,amssymb,amsthm}
\usepackage{enumitem}
\usepackage{xcolor}
\usepackage{hyperref}
\usepackage{caption}

\usepackage{lineno}
\modulolinenumbers[3]

\theoremstyle{plain}

\theoremstyle{definition}
\newtheorem{definition}{Definition}

\theoremstyle{remark}

\newcommand{\C}{\mathbb{C}}
\newcommand{\Q}{\mathbb{Q}}

\numberwithin{equation}{section}

\hypersetup{colorlinks = true, linkcolor = blue, urlcolor = blue, citecolor = blue}

\begin{document}

\title[Unified constructions of the regular Heptagon, Triskaidecagon ...]
      {Unified constructions of the regular Heptagon, Triskaidecagon and Heptadecagon}

\author[H. Ruhland]{Helmut Ruhland}
\address{Santa F\'{e}, La Habana, Cuba}
\email{helmut.ruhland50@web.de}



\begin{abstract}
Constructions of regular heptagon and triskaidecagon by trisection of an angle are well known.
An elegant construction of the heptagon by S. Adlaj shows a 3-fold symmetry related to a Galois group.
Based on the latter construction, in this article one more for the heptagon, two more for the triskaidecagon and three for heptadecagon are presented, all using angle trisection.
\end{abstract}

\date{\today}

\maketitle

\section{Introduction}

In \cite{Gleason} A. Gleason gave constructions of a regular heptagon and triskaidecagon using only square roots and the trisection of angles. Later S. Adlaj \cite{Adlaj} presented a very elegant construction of the heptagon, that shows the action of a cyclic subgroup of order $3$ in the related Galois group on the 3 constructed vertices. \\

In this article I present 3 new unified constructions of the 2 regular polygons by trisection, unified means here \emph{all} based on S. Adlaj's geometric construction. They all use angle trisection.\\

In section \ref{heptagon} I repeat the construction in \cite{Adlaj} and show one more for the heptagon. 
In section \ref{tridecagon} 2 constructions for the trikaidecagon are shown. \\

In section \ref{pentagon} a construction for the pentagon, in section \ref{heptadecagon} $3$ constructions for the heptadecagon are presented.

\section{Constructions of the heptagon \label{heptagon}}

\subsection{The first construction, type I: S. Adlaj's from \cite{Adlaj}}

Let $\epsilon_k = e^{2 \pi i / 3 \, k} = ((-1 + \sqrt{3} \, i) / 2)^k,  \enspace k = 0, 1, 2$ a third root of unity. Define the 3 third roots $\zeta_k,  \enspace k = 0, 1, 2$:
\begin{equation}
   \zeta_k = \epsilon_k \sqrt[3]{\zeta}  \qquad \zeta = \frac{1 - 3 \sqrt{3} \, i}{2 \sqrt{7}}
   \label{zeta_k_7}
\end{equation}
Determining the third root of $\zeta$ is a equivalent to a trisection because $\vert \zeta \vert = 1$. The angle to trisect is $\theta = - \arctan (3 \sqrt{3}) \approx - 79.1066^{\circ}$  \\

Take as radii of the 2 grey, concentric circles in figure \ref{SC}:
\begin{equation}
   R_1, R_2 = \sqrt{(7 \pm \sqrt{21}) / 18}
   \label{R12}
\end{equation}

Form three vertices of the heptagon, the three {\color{red} red}, {\color{green} green}, {\color{blue} blue} parallelograms in figure \ref{SC} realize the following complex additions:
\begin{equation}
   { \color{red}   V_0 = R_1 \epsilon_0 + R_2 \zeta_1 } \qquad
   { \color{green} V_1 = R_1 \epsilon_1 + R_2 \zeta_0 } \qquad 
   { \color{blue}  V_2 = R_1 \epsilon_2 + R_2 \zeta_2 }
   \; \footnote{There is no typing error here, in $V_0, V_1$ the $\epsilon$ and $\zeta$ subscripts are different}
   \label{V123}
\end{equation}

The 3 constructed vertices $V_0, V_1, V_2$ represent the quadratic residues $PR_2 = \{ 1, 2, 4 \}$ modulo $7$. The triangle built by these 3 vertices has 3 different side lengths and so no symmetry. \\
The cyclic permutation $C_3 : \zeta_k \rightarrow \zeta_{k+1}$ of order 3, subscripts modulo 3,
acts on these 3 vertices.  

\subsection{The second construction, Type II \label{Con2_7}}

Now the angle to trisect is $\theta = 0 = 0^{\circ}$. Only a square root is necesary to get $\theta / 3 = 0^{\circ}, 120^{\circ}, 240^{\circ}$. \\

Set $R_1 = 1$ and $R_2$ equals one of the six real roots $r_k, k = 0, \dots, 5$ of the palindromic, sextic  polynomial $P = r^6 + 6 r^5 - 6 r^4 - 29 r^3 - 6 r^2 - 6 r + 1$. Because $P$ is palindromic, with every root $r$, the inverse is a root too. The $6$ roots can be constructed by square roots and solving a cubic, because with $Q = s^3 + 6 s^2 - 9 s - 41$ and $s = r + 1 / r$, $P = Q (s)  \, r^3$. \\
These $6$ roots can also be expressed by the $\zeta_k$ gotten by the trisection \ref{zeta_k_7}: 
\begin{equation}
   \begin{split}
   s_k            & = - 2 + \sqrt{7} \; \left( \zeta_k + \overline{\zeta_k} \right)   \hbox{\hspace{36mm}} \\
   r_{k}, r_{k+3} & = \frac{s_k \pm \sqrt{s_k^2 - 4}}{2} \qquad s_k \; \hbox{a root of } Q \qquad k = 0, 1, 2
   \label{r6_7_tri}
   \end{split}
\end{equation} \\

The vertices of the heptagon are given by \ref{V123}. Hint: a negative $R_2$ in this geometric construction means: a vector in the parallelogram has to be taken negative. \\

In contrary to the previous construction, the triangle built by the 3 constructed vertices $V_0, V_1, V_2$ is now an \emph{isosceles} triangle with a reflection symmetry. \\
$C_3 : \zeta_k \rightarrow \zeta_{k+1}$ acts now on the $6$ radii $r_k$ and so on the $6$ isosceles triangles. Up to scaling and rotation are only 3 different triangles.

\subsection{The third construction, but nothing new, is of type II too}

Now the angle to trisect is $\theta = \pi = 180^{\circ}$. Only a square root is necesary to get $\theta / 3 = 60^{\circ}, 180^{\circ}, 300^{\circ}$.

Set $R_1 = 1$ and $R_2$ equals the \emph{negative} of one of the $6$ real roots $r_k$ of the second construction
in \ref{Con2_7}. \\

Because a negative radius $R_2$ is equivalent to a positive radius and a rotation by $\pi$ of the triangle inscribed in this circle: this construction is equivalent to the previous subsection \ref{Con2_7} . \\

\section{Constructions of the triskaidecagon \label{tridecagon}}

\subsection{The first construction, type I}

Let $\epsilon_k = e^{2 \pi i / 3 \, k} = ((-1 + \sqrt{3} \, i) / 2)^k,  \enspace k = 0, 1, 2$ a third root of unity. Define the 3 third roots $\zeta_k,  \enspace k = 0, 1, 2$:
\begin{equation}
   \zeta_k^{\pm} = \epsilon_k \sqrt[3]{\zeta^{\pm}}  \qquad
      \zeta^{\pm} = \frac{\sqrt{26 \pm 5 \sqrt{13}} - \sqrt{26 \mp 5 \sqrt{13}} \, i}{2 \sqrt{13}}
   \label{zeta_k_13}
\end{equation}
Determining the third root of $\zeta^{\pm}$ is a equivalent to a trisection because $\lvert \zeta^{\pm} \rvert = 1$. The 2 corresponding $2$ angles to trisect are:
\begin{equation}
   \theta^{\pm} = - \arctan \left(\frac{26 \mp 5 \sqrt{13}}{3 \sqrt{3} \sqrt{13}}\right)
   \label{theta12_13}
\end{equation}

Take as radii of the 2 grey, concentric circles and angle to trisect, in figure \ref{SC}:
\begin{equation}
   R_1, R_2 = \sqrt{\frac{\sqrt{13 + \sqrt{13}} \pm \sqrt{5 + \sqrt{13}}}{2 \sqrt{2}}} \quad R_1 \, R_2 = 1 
   \quad \theta^+ \approx - 23.0510^{\circ}
   \label{R12_13}
\end{equation}

The 3 constructed vertices $V_0, V_1, V_2$, see \ref{V123}, represent the now quartic residues $PR_4 = \{ 1, 3, 9 \}$ modulo $13$. Using $- \, \theta^+ / 3$ the vertices represent $4 \, PR_4 =  \{ 4, 10, 12 \}$. The triangle built by these 3 vertices has 3 different side lengths. \\
The cyclic permutation $C_3 : \zeta_k \rightarrow \zeta_{k+1}$ of order 3, subscripts modulo 3,
acts on these 3 vertices. \\

To get the vertices belonging to the $2$ remaing multiplicative cosets of quartic residues $2 \, PR_4 = \{ 2, 5, 6 \}$ and $8 \, PR_4 =  \{ 7, 8, 11 \}$, the following $2$ radii and angle to trisect have to be used:  
\begin{equation}
   R_1, R_2   = \sqrt{\frac{\sqrt{13 - \sqrt{13}} \pm \sqrt{5 - \sqrt{13}}}{2 \sqrt{2}}}
                      \qquad R_1 \, R_2 = 1 \quad \theta^- \approx - 66.9489^{\circ}
   \label{R12_13_other}
\end{equation}

\subsection{The second construction, Type II \label{Con2_13}}

Now the angle to trisect is $\theta = 0 = 0^{\circ}$. Only a square root is necesary to get $\theta / 3 = 0^{\circ}, 120^{\circ}, 240^{\circ}$. \\

Set $R_1 = 1$ and $R_2$ equals one of the 12 real roots $r_k, k = 0, 1, \dots, 11$ of the palindromic, degree $12$  polynomial $P = r^{12} + 12 r^{11} - 12 r^{10} - 274 r^9 -441 r^8 + 441 r^7 + 1275 r^6 + 441 r^5 - 441 r^4 - 274 r^3 - 12 r^2 + 12 r + 1$. Because $P$ is palindromic, with every root $r$, the inverse is a root too. The 12 roots can be constructed by square roots and solving a cubic, because with $Q = s^6 + 12 s^5 - 18 s^4 - 334 s^3 - 384 s^2 + 1323 s + 2131$ and $s = r + 1 / r$, $P = Q (s)  \, r^6$. The sectic $Q$ factorizes in $\Q (\sqrt{13})$ as product of $R = s^3 + 3 \, (2 + \sqrt{13}) \, s^2 + 21 \, (3 + 1 \sqrt{13}) / 2 \, s + (15 + 107 \sqrt{13}) / 2)$ and its $\Q (\sqrt{13})$-conjugate $\overline{R}$. \\
These $12$ roots can be expressed by the $\zeta_k^{\pm}$ gotten by the trisection \ref{zeta_k_13}: 
\begin{equation}
   \begin{split}
   s_k^{\pm}        & = - (2 \pm \sqrt{13}) 
                        + \sqrt{13 \pm \sqrt{13}} \; \left( \zeta_k^{\pm} + \overline{\zeta_k^{\pm}} \right)
                        \hbox{\hspace{16mm}} \\
   r_{k}, r_{k+3}   & = \frac{s_k^+ \pm \sqrt{(s_k^+)^2 - 4}}{2} \qquad s_k^+ \; \hbox{a root of } R
                        \qquad k = 0, 1, 2 \\  
   r_{k+6}, r_{k+9} & = \frac{s_k^- \pm \sqrt{(s_k^-)^2 - 4}}{2} \qquad s_k^- \; \hbox{a root of } \overline{R}
                        \qquad k = 0, 1, 2  
   \label{r12_13}
   \end{split}
\end{equation}
The conjugation on the rhs of $s_k^{\pm}$ is now a $\C$-conjugation. \\

The vertices of the trikaidecagon are also given by \ref{V123}. Hint: a negative $R_2$ in this geometric construction means: a vector in the parallelogram has to be taken negative. \\

In contrary to the previous construction, the triangle built by the 3 constructed vertices $V_0, V_1, V_2$ is now an \emph{isosceles} triangle with a reflection symmetry. \\
$C_3 : \zeta_k \rightarrow \zeta_{k+1}$ acts now on the $12$ radii $r_k$ and so on the $12$ isosceles triangles. Up to scaling and rotation are only 6 different triangles.

\subsection{The third construction, but nothing new, is of type II too}

Now the angle to trisect is $\theta = \pi = 180^{\circ}$. Only a square root is necesary to get $\theta / 3 = 60^{\circ}, 180^{\circ}, 300^{\circ}$.

Set $R_1 = 1$ and $R_2$ equals the \emph{negative} of one of the $12$ real roots $r_k$ of the second construction
in \ref{Con2_13}. \\

Because a negative radius $R_2$ is equivalent to a positive radius and a rotation of the triangle inscribed in this circle: this construction is equivalent to the previous subsection \ref{Con2_13} . \\

\section{Constructions of the pentagon \label{pentagon}}

\subsection{The first construction, n o \,  type I}

A type I construction does not exist. Every triangle built by vertices in a regular pentagon is an isoceles triangle. A triangle with 3 different side lengths for type I does not exist.

\subsection{The second construction, Type II \label{Con2_5}}

Now the angle to trisect is $\theta = 0 = 0^{\circ}$. Only a square root is necesary to get $\theta / 3 = 0^{\circ}, 120^{\circ}, 240^{\circ}$. \\

Set $R_1 = 1$ and $R_2$ equals one of the four real roots $r_k, k = 0, \dots, 3$ of the palindromic, quartic  polynomial $P = r^4 - r^3 - 9 r^2 - r + 1$. Because $P$ is palindromic, with every root $r$, the inverse is a root too. The $4$ roots can be constructed by square roots, because with $Q = s^2 - s - 11$ and $s = r + 1 / r$, $P = Q (s)  \, r^2$. No trisection is needed. \\
\begin{equation}
   \begin{split}
   s_0, s_1         & = (1 \pm 3\, \sqrt{5}) / 2  \qquad \hbox{the 2 roots of } Q \\
   r_{k}, r_{k+2}   & = \frac{s_k \pm \sqrt{s_k^2 - 4}}{2}  \qquad k = 0, 1 \hbox{\hspace{16mm}}
   \label{r4_5}
   \end{split}
\end{equation}

The vertices of the pentagon are given by \ref{V123}. Hint: a negative $R_2$ in this geometric construction means: a vector in the parallelogram has to be taken negative. \\

The triangle built by the 3 constructed vertices $V_0, V_1, V_2$ is an \emph{isosceles} triangle with a reflection symmetry. \\

\subsection{The third construction, but nothing new, is of type II too}

Now the angle to trisect is $\theta = \pi = 180^{\circ}$. Only a square root is necesary to get $\theta / 3 = 60^{\circ}, 180^{\circ}, 300^{\circ}$.

Set $R_1 = 1$ and $R_2$ equals the \emph{negative} of one of the $4$ real roots $r_k$ of the second construction in \ref{Con2_5}. \\

Because a negative radius $R_2$ is equivalent to a positive radius and a rotation by $\pi$ of the triangle inscribed in this circle: this construction is equivalent to the previous subsection \ref{Con2_5} . \\

\section{Constructions of the heptadecagon \label{heptadecagon}}

In 1796 Carl Friedrich Gauss showed the algebraic part of a construction of the heptagon, see \cite{Gauss},
using only square roots. The incidental geometric part using a compass and unmarked straightedge was
not given for almost a century.
With Galois theory it can be shown: if $2^{2^n} + 1$ is a prime $p$, then the regular $p$-gon can be constructed using a compass and unmarked straightedge. Primes of this form are called Fermat primes. Algebraically spoken: the field extension of $\Q$ by the $p$-th unit roots has degree $2^{2^n}$
and a cyclic group as Galois group. It can be constructed by $2^n$ subsequent extensions of degree $2$
by an adjunction of square roots. For $n = 2$ we get $p = 17$ and so the heptadecagon or $17$-gon
can be constructed by compass and unmarked straightedge or subsequent field extensions of square roots. \\

\begin{definition}
Let $1 \le a \le b \le c \le n-2$. A triangle built by $3$ vertices of a $n$-gon is denoted by ($a, b, c$ the triangle "side lengths"): \\
- a list $[ a, b, c ]$ if $a + b - c = 0$, a list because $c$ is a distinguished side length \\
- a set $\{ a, b, c \}$ if $a + b + c = 0$, a set because the condition is symmetric in $a, b, c$ \\
The "side lengths" are given as the angle under which the side is seen from the midpoint of the $n$-gon,
the midpoint of the circumscribed circle, in multiples of $2 \pi / n$. In the case of a $[ \dots ]$ triangle the midpoint lies outside of the triangle. This denotation of triangles is invariant under scaling, rotation and reflection.
\end{definition}
 
\subsection{The first construction, type I}

The construction by the method of S. Adlaj uses trisection of an angle, which is equivalent to
a field extension of degree $3$. $\Q$ extended by the 17-th unit roots has as degree a power of $2$.
So any construction using angle trisection, equivalent to a degree $3$ extension, seems not to be
useful to get a power of $2$ field extension. For the current construction it turns out, that the polynomial
of degree $3 \cdot 32 = 96$, the roots defining a trisected angle $\theta / 3$, factorizes over $\Q$.
We get $3$ polynomials of degree $32$ and the $3$ angles $\theta / 3, \theta / 3 + 120^{\circ}, \theta / 3 + 240^{\circ}$ as result of the trisection are each defined by one of the $3$ polynomials. \\

Let $\epsilon_k = e^{2 \pi i / 3 \, k} = ((-1 + \sqrt{3} \, i) / 2)^k,  \enspace k = 0, 1, 2$ a third root of unity. Define the $3$ $\zeta_k,  \enspace k = 0, 1, 2$ with $\lvert \zeta_k \rvert = 1$, representing
the $3$ trisected angles $\theta / 3, \theta / 3 + 120^{\circ}, \theta / 3 + 240^{\circ}$:
\begin{equation}
  \zeta_k = \epsilon_k \zeta    \qquad \zeta \; \hbox{a root of the polynomial } Z, \lvert \zeta \rvert = 1 
   \label{zeta_k_17}
\end{equation}

Take as radii of the 2 grey, concentric circles in figure \ref{SC}:
\begin{equation}
   R_1 = 1 / \sqrt{r}, R_2 = \sqrt{r} \quad R_1 \, R_2 = 1 
   \qquad r \; \hbox{a root of the polynomial } R
   \label{R12_13}
\end{equation}

The $2$ polynomials $Z$ and $R$ used above are defined for each of the $2$ cases in the following $2$ subsections. To construct the $3$ vertices $V_0, V_1, V_2$, see \ref{V123}.

\subsubsection{The type Ia construction}

The following polynomial $Z$ defines the $\zeta_k$ in \ref{zeta_k_17}. $Z$ is palindromic. All
roots of $Z$ are complex with an absolute value of $1$. The roots of $Z_R$ are the real parts of the roots of $Z$. Because for a root of $Z$ the complex conjugate with the same real part is also a root, the degree of $Z_R$ is only half of the degree of $Z$:
\begin{equation}
   \begin{split}
   Z = 	\enspace & 1429 \zeta^{32} + 4724 \zeta^{30} + 5046 \zeta^{28} + 952 \zeta^{26} - 1753 \zeta^{24}
         - 3012 \zeta^{22} - 4328 \zeta^{20} \\
       - \; & 1024\zeta^{18} + 2493 \zeta^{16} - 1024 \zeta^{14} - 4328 \zeta^{12} - 3012 \zeta^{10}
         - 1753 \zeta^8 + 952 \zeta^6 \\
       + \; & 5046 \zeta^4 + 4724 \zeta^2 + 1429 \\
   Z' = 	\enspace & 1429 \zeta^{32} - 3376 \zeta^{30} + \dots  \qquad
   Z'' = 	1429 \zeta^{32} - 1348 \zeta^{30} + \dots \\
   Z_R = \enspace & 153812992 \zeta_R^{16} - 618758144 \zeta_R^{14} + 1013051392 \zeta_R^{12}
         - 865433600 \zeta_R^{10} \\
       + \; & 410237440 \zeta_R^8 - 105536768 \zeta_R^6 + 13271872 \zeta_R^4 - 643112 \zeta_R^2 + 4489
   \label{IaM_17}
   \end{split}
\end{equation}
As already mentioned above, the $3$ roots corresponding to the $3$ trisected angles $\theta / 3, \theta / 3 + 120^{\circ}, \theta / 3 + 240^{\circ}$ are not roots of only $Z$, but roots of $Z, Z', Z''$.

\begin{equation}
   \begin{split}
   R = \enspace & 1429 r^{32} - 66724 r^{30} + 1243602 r^{28} - 12462308 r^{26} + 75668174 r^{24} \\
- \; & 295963182 r^{22} + 771231283 r^{20} - 1363022242 r^{18} + 1646851473 r^{16} \\
- \; & 1363022242 r^{14} + 771231283 r^{12} - 295963182 r^{10} + 75668174 r^8 \\
- \; & 12462308 r^6 + 1243602 r^4 - 66724 r^2 + 1429
   \label{IaR_17}
   \end{split}
\end{equation}

Up to scaling, rotation and reflection are $8$ different triangles: $[ 1, 2, 3 ], [ 1, 5, 6 ], [ 1, 7, 8 ],$
$[ 2, 4, 6 ], [ 2, 5, 7 ], [ 3, 4, 7 ]$ with the midpoint outside of the triangles and 
$\{ 3, 6, 8 \}, \{ 4, 5, 8 \}$.

\subsubsection{The type Ib construction}

\begin{equation}
   \begin{split}
   Z = \enspace & 2347 \zeta^{32} - 214 \zeta^{30} + 243 \zeta^{28} + 388 \zeta^{26} - 712 \zeta^{24}
       - 972 \zeta^{22} + 1249 \zeta^{20} \\
      + \; & 1316 \zeta^{18} - 729 \zeta^{16} + 1316 \zeta^{14} + 1249 \zeta^{12} - 972 \zeta^{10}
        - 712 \zeta^8 + 388 \zeta^6 \\
      +  \; & 243 \zeta^4 - 214 \zeta^2 + 234793650944 \\
	 Z_R = \enspace & 93650944 \zeta_R^{16} - 297205760 \zeta_R^{14} + 358506496 \zeta_R^{12} - 203632640 \zeta_R^{10} \\
         + \; & 54351616 \zeta_R^8 - 5904128 \zeta_R^6 + 242464 \zeta_R^4 - 2432 \zeta_R^2 + 1
   \label{IbM_17}
   \end{split}
\end{equation}

\begin{equation}
   \begin{split}
   R = \enspace &	2347 r^{32} - 119050 r^{30} + 2192967 r^{28} - 21082328 r^{26} + 122159825*r^{24} \\
- \; & 458344224 r^{22} + 1157048782 r^{20} - 2004900694 r^{18} + 2406196287 r^{16} \\
- \; & 2004900694 r^{14} + 1157048782 r^{12} - 458344224 r^{10} + 122159825 r^8 \\
- \; & 21082328 r^6 + 2192967 r^4 - 119050 r^2 + 2347 
   \label{IbR_17}
   \end{split}
\end{equation}

Up to scaling, rotation and reflection are $8$ different triangles: $[ 1, 3, 4 ], [ 1, 4, 5 ], [ 1, 6, 7 ],$
$[ 2, 3, 5 ], [ 2, 6, 8 ], [ 3, 5, 8 ]$ with the midpoint outside of the triangles and 
$\{ 2, 7, 8 \}, \{ 4, 6, 7 \}$.

\subsection{The second construction, Type II \label{Con2_17}}

Now the angle to trisect is $\theta = 0 = 0^{\circ}$. Only a square root is necesary to get $\theta / 3 = 0^{\circ}, 120^{\circ}, 240^{\circ}$. \\

Set $R_1 = 1$ and $R_2$ equals one of the 16 real roots $r_k, k = 0, 1, \dots, 15$ of the palindromic, degree $16$  polynomial $P = r^{16} - r^{15} - 135 r^{14} - 545 r^{13} + 545 r^{12} + 5643 r^{11}
+ 6733 r^{10} - 6733 r^9 - 17577 r^8 - 6733 r^7 + 6733 r^6 + 5643 r^5 + 545 r^4 - 545 r^3 - 135 r^2 - r + 1$. Because $P$ is palindromic, with every root $r$, the inverse is a root too. The 16 roots can be constructed by square roots, because with $Q = s^8 - s^7 - 143 s^6 - 538 s^5 + 1375 s^4 + 8354 s^3
+ 3322 s^2 - 26380 s - 29681$ and $s = r + 1 / r$, $P = Q (s)  \, r^8$. The octic $Q$ has the cyclic group of order $8$ as Galois group, this group is solvable and the $8$ roots $s_0, \dots, s_7$ can be expressed using only square roots. \\

The $4$ symbols $\alpha, \beta, \gamma, \delta$ represent $\pm1$ signs for the 4 different roots, but not
all  of these $2^4 = 16$ combinations are allowed, because $Q$ has only $8$ roots. The $3$ signs $\alpha, \beta, \gamma$ are free, $\delta = \alpha \, \beta$, so we get $2^3 = 8$ roots:
\begin{equation}
   \begin{split}
   & \quad \enspace A = \alpha \, \sqrt{2 \, (85 - \delta \, 19 \, \sqrt{17})}  \qquad
                    B = \beta  \, \sqrt{2 \, (85 - \delta \, 13 \, \sqrt{17})} \\
   s_0, \dots, s_7 & = \frac{1 - \delta \, 3 \, \sqrt{17} + 3 \, A
                       - \gamma \, \sqrt{3 \, (204 - \delta \, 36 \, \sqrt{17}
                       + 2 \, (2 - \delta \, 6 \, \sqrt{17}) \, A + 16 \, B )}}{8} \\
   r_{k}, r_{k+8}  & = \frac{s_k \pm \sqrt{s_k^2 - 4}}{2} \qquad s_k \; \hbox{a root of } Q
                        \qquad k = 0, \dots, 7  
   \label{s8_17}
   \end{split}
\end{equation}

The vertices of the heptadecagon are also given by \ref{V123}. Hint: a negative $R_2$ in this geometric construction means: a vector in the parallelogram has to be taken negative. \\

In contrary to the previous construction, the triangle built by the 3 constructed vertices $V_0, V_1, V_2$ is now an \emph{isosceles} triangle with a reflection symmetry. \\
Up to scaling, rotation and reflection are $8$ different triangles: $[ 1, 1, 2 ], [ 2, 2, 4 ], [ 3, 3, 6 ],$ $[ 4, 4, 8 ]$ with the midpoint outside the triangles and 
$\{ 1, 8, 8 \}, \{ 3, 7, 7 \}, \{ 5, 6, 6, \}, \{ 5, 5, 7 \}$.

\subsection{The third construction, but nothing new, is of type II too}

Now the angle to trisect is $\theta = \pi = 180^{\circ}$. Only a square root is necesary to get $\theta / 3 = 60^{\circ}, 180^{\circ}, 300^{\circ}$.

Set $R_1 = 1$ and $R_2$ equals the \emph{negative} of one of the $16$ real roots $r_k$ of the second construction in \ref{Con2_17}. \\

Because a negative radius $R_2$ is equivalent to a positive radius and a rotation of the triangle inscribed in this circle: this construction is equivalent to the previous subsection \ref{Con2_17} . \\

\section{Open questions}

Are the constructions of type I, type II or both also possible for other p-gons with $p$ a prime of the form $6 n + 1$? The next primes $p$ to investigate would be $p = 19, 31, \dots$.

\newpage
\noindent \textbf{\large Appendices}

\appendix

\section{The geometric construction in \cite{Adlaj}}

\begin{center}
  \includegraphics[width=0.75\textwidth]{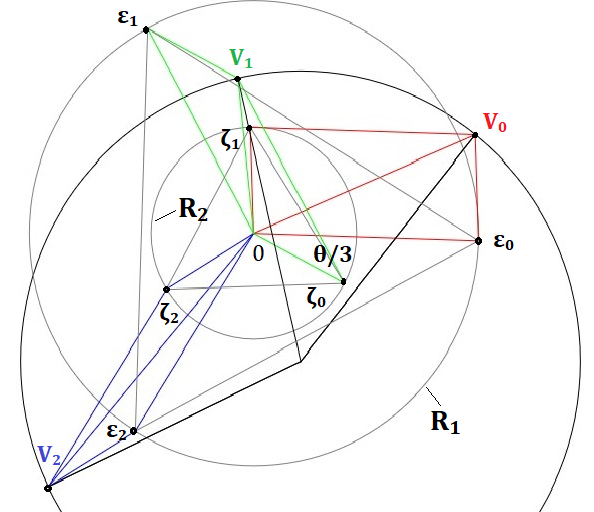}
  \captionof{figure}{A scan of the heptagon construction in \cite{Adlaj} with the denomination
                     of vertices added by the author of this article. The vertices $\epsilon_k, \zeta_k$
                     do not have unit distance to $0$, instead they should have been denominated $\epsilon_k R_1, \zeta_k R_2$}
  \label{SC}
\end{center}

\bibliographystyle{amsplain}

\end{document}